\documentclass[runningheads]{llncs}

\usepackage{amssymb,amsmath}
\usepackage{graphicx}

\usepackage{url}
\urldef{\mailsa}\path|{jchifman, petrovic}@ms.uky.edu|

\newtheorem{alg}[theorem]{Algorithm}

\newcommand{\ul} \underline

\newcommand{\mtx}{\left[ \begin{matrix}}
\newcommand{\mtxend}{\end{matrix}\right]}

\newcommand{\rk}{\operatorname{rank}}
\newcommand{\cols}{\operatorname{cols}}

\begin{document}

\mainmatter

\title{Toric ideals of phylogenetic invariants for the general group-based model on  claw trees $K_{1,n}$ }
\titlerunning{Phylogenetic ideals for claw trees}
\author{Julia Chifman, Sonja Petrovi\'{c}}
\authorrunning{J. Chifman, S. Petrovi\'{c}}
\institute{Department of Mathematics, University of Kentucky, Lexington, KY 40506, USA
\mailsa}

\toctitle{Phylogenetic ideals for claw trees}
\tocauthor{J. Chifman, S. Petrovi\'{c}}
\maketitle

\begin{abstract}
We address the problem of studying the toric ideals of phylogenetic invariants for a general group-based model on an arbitrary claw tree. We focus on the group
$\mathbb Z_2$ and choose a natural recursive approach that extends to other groups.  The study of the lattice associated with each phylogenetic ideal produces a
list of circuits that generate the corresponding lattice basis ideal.  In addition, we describe {\em explicitly} a quadratic lexicographic Gr\"{o}bner basis of
the toric ideal of invariants for the claw tree on an arbitrary number of leaves. Combined with a result of Sturmfels and Sullivant, this implies that the
phylogenetic ideal of {\em every} tree for the group $\mathbb Z_2$ has a quadratic Gr\"{o}bner basis. Hence, the coordinate ring of the toric variety is a Koszul
algebra.

\end{abstract}

\subsubsection*{Acknowledgment.} The authors would like to thank Uwe Nagel for introducing us to the field of phylogenetic
 algebraic geometry and for his continuous support, motivation and guidance.

\section{Introduction}

Phylogenetics is concerned with determining genetic relationship between species based on their DNA sequences.  First, the various DNA sequences are aligned,
that is, a correspondence is established that accounts for their differences.  Assuming that all DNA sites evolve identically and independently, the focus is on
one site at a time.  The data then consists of observed pattern frequencies in aligned sequences. This observed data are used to estimate the true  joint
probabilities of the observations and, most importantly, to reconstruct the ancestral relationship among the species. The relationship can be represented by a
phylogenetic tree.

A {\emph{phylogenetic tree}} $T$ is a simple, connected, acyclic graph equipped with some statistical information. Namely, each node of $T$ is a random variable
with $k$ possible states chosen from the state space $S$. Edges of $T$ are labeled by transition probability matrices that reflect probabilities of changes of
the states from a node to its child. These probabilities of mutation  are the parameters for the statistical model of evolution, which is described in terms of a
discrete-state continuous-time Markov process on the tree.
 Since the goal is to reconstruct the tree, the interior nodes are hidden.  The relationship between
the random variables is encoded by the structure of the tree. At each of the $n$ leaves, we can observe any of the $k$  states; thus there are $k^n$ possible
observations.   Let $p_\sigma$ be the joint probability of making a particular observation $\sigma\subset S^n$ at the leaves.  Then $p_\sigma$ is a polynomial in
the model parameters.

A {\emph{phylogenetic invariant}} of the model is a polynomial in the leaf probabilities which vanishes for every choice of model parameters. The set of these
polynomials forms a prime ideal in the polynomial ring over the unknowns $p_{\sigma}$. The objective is to compute this ideal explicitly.    Thus we consider a
polynomial map  $ \phi:\mathbb C^N \to \mathbb C^{k^n}$, where $N$ is the total number of model parameters. The map depends only on the tree $T$ and the number
of states $k$; its coordinate functions are the $k^n$ polynomials $p_\sigma$. The map $\phi$ induces a parametrization of an algebraic variety.  The study of
these algebraic varieties for various statistical models is a central theme in the field of algebraic statistics (\cite{StSu}). Phylogenetic invariants are a
powerful tool for tree reconstruction (\cite{AlRh2}, \cite{CaFS}, \cite{ErYa}).

There is a specific class of models for which the ideal of invariants is particularly nice. Let $M_e$ be the $k\times k$ transition probability matrix for edge
$e$ of $T$. In the general Markov model, each matrix entry is an independent model parameter.   A group-based model is one in which the matrices $M_e$ are
pairwise distinct, but it is required that certain entries coincide. For these models, transition matrices are diagonalizable by the Fourier transform of an
abelian group.  The key idea behind this linear change of coordinates is to label the states (for example, $A$,$C$,$G$, and $T$) by a finite abelian group (for
example, $\mathbb Z_2\times \mathbb Z_2$) in such a way that transition from one state to another depends only on the difference of the group elements. Examples
of group-based models include the Jukes-Cantor and Kimura's one-parameter models used in computational biology.

Sturmfels and Sullivant in \cite{StSu} reduce the computation of ideals of phylogenetic invariants of group-based models on an arbitrary tree to the case of claw
trees $T_n:=K_{1,n}$, the complete bipartite graph from one node (the root) to $n$ nodes (the leaves).  The main result of \cite{StSu} gives a way of
constructing the ideal of phylogenetic invariants for any tree {\em if} the ideal for the claw tree is known.  However, in general, it is an open problem to
compute the phylogenetic invariants for a claw tree. We consider the ideal for a general group-based model for the group $\mathbb Z_2$. Let $q_\sigma$ be the
image of $p_\sigma$ under the Fourier transform.  Assuming the identity labeling function and adopting the notation of \cite{StSu}, the ideal of phylogenetic
invariants for the tree $T_n$ is the kernel of the following homomorphism between polynomial rings:
\begin{align*}
\varphi_n: \mathbb C[q_{g_1,\dots,g_n}:g_1,\dots,g_n\in G] & \to \mathbb C[a_g^{(i)}: g\in G,i=1,\dots,n+1]   \\ q_{g_1,\dots,g_n}  & \mapsto
a_{g_1}^{(1)}a_{g_2}^{(2)}\dots a_{g_n}^{(n)}a_{g_1+g_2+\dots +g_n}^{(n+1)} , \tag{*}\label{phidefn}
\end{align*}
where $G$ is a finite group with $k$ elements, each corresponding to a state.  The coordinate $q_{g_1,\dots,g_n}$ corresponds to observing the element $g_1$ at
the first leaf of T, $g_2$ at the second, and so on.  The phylogenetic invariants form a {\emph{toric ideal}} in the Fourier coordinates $q_\sigma$, which can be
computed from the corresponding lattice basis ideal by saturation. The main result of this paper is a complete description of the lattice basis ideal and a
quadratic Gr\"{o}bner basis of the ideal of invariants for the group  $\mathbb Z_2$ on $T_n$ for any number of leaves $n$.

Our paper is organized as follows.
 In section 2 we lay the foundation for our recursive approach.  The ideal of the two-leaf claw tree is trivial, so we begin
with the case when the number of leaves is three.    Sections 3 and 4 address the problem of describing the lattices corresponding to the toric ideals. We
provide a nice lattice basis consisting of circuits.  The corresponding lattice basis ideal is generated by circuits of degree two and thus in particular
satisfies the Sturmfels-Sullivant conjecture.

The ideal of phylogenetic invariants is the saturation of the lattice basis ideal.  However, we do not use any of the standard algorithms to compute saturation
(e.g. \cite{HeMa}, \cite{St}). Instead, our recursive construction of the lattice basis ideals can be extended to give the full ideal of invariants, which we
describe in the final section. The recursive description of these ideals depends only on the number of leaves of the claw tree and it does not require
saturation.   Finally, and possibly somewhat surprisingly, we show that the ideal of invariants for every claw tree admits a quadratic Gr\"{o}bner basis with
respect to a lexicographic term order. We describe it {\em explicitly}.

Combined with the main result of Sturmfels and Sullivant in \cite{StSu}, this implies that the phylogenetic ideal of {\em every} tree for the group $\mathbb Z_2$
has a quadratic Gr\"{o}bner basis. Hence, the coordinate ring of the toric variety is a Koszul algebra.  In addition, the ideals for every tree can be computed
explicitly. These ideals are particularly nice as  they satisfy the conjecture in \cite{StSu} which proposes that the order of the group gives an upper bound for
the degrees of minimal generators of the ideal of invariants.  The case of $\mathbb Z_2$ has been solved in \cite{StSu} using a technique that does not
generalize.  We hope to extend our recursive approach and obtain the result for an arbitrary abelian group.

For a detailed background on phylogenetic trees, invariants, group-based models, Fourier coordinates, labeling functions and more, the reader should refer to
\cite{AlRh}, \cite{ErSt}, \cite{PaSt}, \cite{StSu}.

\section{Matrix representation}

\paragraph{}
Fix a claw tree $T_n$ on $n$ leaves and a finite abelian group $G$ of order $k$.  Soon we will specialize to the case $k=2$. We want to compute the ideal of
phylogenetic invariants for the general group-based model on $T_n$.  After the Fourier transform, the ideal of invariants (in Fourier coordinates) is given by
$I_n=\ker\varphi_n$ , where $\varphi_n$ is a map between polynomial rings in $k^n$ and $k(n+1)$ variables, respectively, defined by (\ref{phidefn}).    In order
to compute the toric ideal $I_n$, we first compute the lattice basis ideal $I_{L_n}\subset I_n$ corresponding to $\varphi_n$ as follows.  Fixing an order on the
monomials of the two polynomial rings, the linear map $\varphi$ can be represented by a matrix $B_{n,k}$ that describes the action of $\varphi$ on the variables.
Then the lattice ${L_n}=\ker(B_{n,k})\subset\mathbb Z^{k^n}$ determines the ideal $I_{L_n}$. It is generated by elements of the form $ (\prod q_{g_1,\dots,g_n})
^ {v^+} -  (\prod q_{g_1,\dots,g_n})^{v^-} $ where $v=v^+ -v^-\in {L_n}$. We will give an explicit description of this basis and, equivalently, the ideals
$I_{L_n}$.

Hereafter assume that $G=\mathbb Z_2$. For simplicity, let us say that $B_n:=B_{n,2}$.

To create the matrix $B_n$, first order the two bases as follows.  Order the $a_g^{(i)}$ by varying the upper index $(i)$ first and then the group element $g$:
$a_0^{(1)}$, $a_0^{(2)}$, $\dots$, $a_0^{(n+1)}$, $a_1^{(1)}$, $\dots$, $a_{1}^{(n+1)}$. Then, order the $q_{g_1,\dots,g_n}$ by ordering the indices with respect
to binary counting: $$q_{0\dots 00}> q_{0\dots 01}> \dots> q_{1\dots 10}>q_{1\dots 1}.$$ That is, $q_{g_1\dots g_n}>q_{h_1\dots h_n}$ if and only if $(g_1\dots
g_n)_2 < (h_1\dots h_n)_2$, where $$(g_1\dots g_n)_2:=g_12^{n-1}+g_22^{n-2}+\dots +g_n2^0$$ represents the binary number $g_1\dots g_n$.\\

Next, index the rows of $B_{n}$ by $a_g^{(i)}$ and its columns by $q_{g_1,\dots,g_n}$.  Finally, put  $1$ in the entry of $B_{n}$ in the row indexed by
$a_g^{(i)}$ and column indexed by $q_{g_1,\dots,g_n}$ if  $a_g^{(i)}$ divides the image of  $q_{g_1,\dots,g_n}$, and $0$ otherwise.

\begin{example}
Let $n=2$. Then we order the $q_{ij}$ variables according to binary counting: $q_{00}$, $q_{01}$, $q_{10}$, $q_{11}$, so that
\begin{align*}
\varphi: \mathbb C[q_{00}, q_{01}, q_{10}, q_{11}] &\to \mathbb C[a_0^{(1)},a_0^{(2)},a_0^{(3)},a_1^{(1)},a_1^{(2)},a_1^{(3)}]\\ q_{00} &\mapsto
a_{0}^{(1)}a_{0}^{(2)}a_{0+0}^{(3)}\\ q_{01} &\mapsto a_{0}^{(1)}a_{1}^{(2)}a_{0+1}^{(3)}\\ q_{10} &\mapsto a_{1}^{(1)}a_{0}^{(2)}a_{1+0}^{(3)}\\ q_{11} &\mapsto
a_{1}^{(1)}a_{1}^{(2)}a_{1+1}^{(3)}.
\end{align*}
Now we put the $a_{i}^{(j)}$ variables in order: $a_0^{(1)}$, $a_0^{(2)}$, $a_0^{(3)}$, $a_1^{(1)}$, $a_1^{(2)}$, $a_1^{(3)}$ .  Thus
$$B_2 = \left[
\begin{matrix}
1 &1 &0 &0\\
1 &0 &1 &0\\
1 &0 &0 &1\\
0 &0 &1 &1\\
0 &1 &0 &1\\
0 &1 &1 &0
\end{matrix}\right].
$$
\end{example}

The tree $T_{n-1}$ can be considered as a subtree of $T_n$ by ignoring, for example,  the leftmost leaf of $T$. As a consequence, a natural question arises: how
does $B_{n}$ relate to $B_{n-1}$?
\begin{remark}\label{obs:1}
The matrix $B_{n-1}$ for the subtree of $T_n$ with the leaf $(1)$ removed can be obtained as a submatrix of $B_n$ for the tree $T_n$ by deleting rows  $1$ and
$(n+1)+1$ and taking only the first $2^{n-1}$ columnns.\\ Divide the $n$-leaf matrix $B_n$ into a $2\times 2$ block matrix with blocks of size $(n+1)\times
2^{n-1}$: $$B_n=\left[ \begin{matrix} B_{11} & B_{12}\\ B_{21} & B_{22}
\end{matrix} \right].$$
Then, grouping together $B_{11},B_{21}$ without the first row of each $B_{i1}$, we obtain the matrix $B_{n-1}$.  This is true because  rows $1$ and $(n+1)+1$
represent the variables $a^{(1)}_{g}$ for $g\in G$ associated with the leaf $(1)$ of $T_n$. Note that the entries in row $a_g^{(n+1)}$ remain undisturbed as the
omitted rows are indexed by the identity of the group.
\end{remark}

\begin{example}
The matrix $B_{2}$ is  equal to the submatrix of $B_{3}$ formed by rows $2$,$3$,$4$,$6$,$7$,$8$, and first $4$ columns.
\end{example}

\begin{remark}\label{obs:2}
Fix any observation $\sigma = g_1,\dots,g_n$ on the leaves. Clearly, at any given leaf $j\in\{1,\dots,n\}$, we observe exactly one group element, $g_j$. Since
the matrix entry     $ b_{a_{g_j}^{(j)},q_\sigma}$     in the row indexed by $a_{g_j}^{(j)}$ and column indexed by $q_\sigma$ is $1$ exactly when $a_{g_j}^{(j)}$
divides the image of $q_\sigma$, one has that $$ \sum_{{g_j}\in G} b_{a_{g_j}^{(j)},q_\sigma} = 1 $$ for a fixed leaf $(j)$ and fixed observation $\sigma$. Note
that the formula also holds if $j=n+1$ by definition of $a_{g_{n+1}}^{(n+1)}=a_{g_1+\dots +g_n}^{(n+1)}$. In particular, the rows indexed by $a_{g_j}^{(j)}$ for
a fixed $j$ sum up to the row of ones.
\end{remark}

\section{Number of lattice basis elements}

\paragraph{} We compute the dimension of the kernel of $B_n$ by induction on $n$. We proceed in two steps.

\begin{lemma} [Lower bound] \label{lm:no1}
$$\rk(B_n)\geq \rk(B_{n-1})+1.$$
\end{lemma}
\begin{proof}
First note that $\rk(B_n)\geq \rk(B_{n-1})$ since $B_{n-1}$ is a submatrix of the first $2^{n-1}$ columns of $B_n$. In the block  $\left[ \begin{matrix} B_{11},
B_{12}\end{matrix}\right]^T$, the row indexed by $a_1^{(1)}$ is zero, while in the  block $\left[\begin{matrix} B_{21},B_{22}\end{matrix}\right]^T$, the row
indexed by  $a^{(1)}_{1}$ is $1$. Choosing one column from $\left[\begin{matrix} B_{21},B_{22}\end{matrix}\right]^T$ provides a vector independent of the first
$2^{n-1}$ columns.  The rank must therefore increase by at least $1$. \hfill$\square$ \end{proof}

\begin{lemma}[Upper bound]\label{lm:no2}
$$\rk(B_n)\leq n+2 .$$
\end{lemma}
\begin{proof}  $B_n$ has $2(n+1)$ rows.  Remark \ref{obs:2} provides $n$ independent relations among the rows of our matrix:  varying $j$ from $1$ to $n+1$, we
obtain that the sum of the rows $j$ and $n+1+j$ is $1$ for each $j=1,\dots,n+1$.   Thus the upper bound is immediate. \hfill$\square$ \end{proof}

We are ready for the main result of the section.
\begin{proposition}[Cardinality of lattice basis]\label{prop:2}

Let $n\geq 2$. Then there are $ 2^n-2(n+1)+n $ elements in the basis of the lattice ${L_n}$ corresponding to $T_n$. That is, $$\dim \ker (B_n) = 2^n-2(n+1)+n. $$
\end{proposition}
\begin{proof}  We show $ \rk(B_n)=2(n+1)-n$. It can be checked directly that $B_2$ has full rank.  Assume that the claim is true for $n-1$.  Then by Lemmae
(\ref{lm:no1}) and (\ref{lm:no2}), $$2(n+1)-n\geq \rk(B_n)\geq \rk(B_{n-1})+1=2n-(n-1)+1 ,$$ where the last equality is provided by the induction hypothesis.
The claim follows since the left- and the right-hand sides agree. \hfill$\square$ \end{proof}

\section{Lattice basis}

\paragraph{} In this section we describe a basis of the kernel of $B_n:=B_{n,2}$, in which the binomials corresponding to the basis elements satisfy the
conjecture on the degrees of the generators of the phylogenetic ideal.  In particular, since the ideal is generated by squarefree binomials and contains no
linear forms, these elements are actually circuits.  By Proposition \ref{prop:2}, we need to find $2^n-(n+2)$ linearly independent vectors in the lattice. The
matrix of the tree with $n=2$ leaves has a trivial kernel, so we begin with the tree on $n=3$ leaves. The dimension of the kernel is $3$ and the lattice basis is
given by the rows of the following matrix:
$$ \left[\begin{matrix} 0&0&1&-1&-1&1&0&0\\ 0&1&0&-1&-1&0&1&0\\1&0&0&-1&-1&0&0&1\end{matrix}\right].$$ In order to
study the kernels of $B_n$ for any $n$, it is useful to have an algorithmic way of constructing the matrices.

\begin{alg}\rm [The construction of $B_n$]\\
{\bf Input:} the number of leaves $n$ of the claw tree $T_n$.\\ {\bf Output:} $B_n\in \mathbb Z^{2(n+1)\times 2^n}$.\\ Initialize $B_n$ to the zero matrix.
\\ Construct the first $n$ rows: \\  \phantom{x}for $k$ from $1$ to $n$ do:\\ \phantom{xx}for $c$ from $0$ to $2^k-1$ with $c\equiv 0 \mod 2$ do:\\
\phantom{xxx}for $j$ from $c2^{n-k}+1$ to $(c+1)2^{n-k}$ do: \phantom{xxx} $b_{k,j}:=1. $ \\ Construct row $n+1$:
\\\phantom{x}if $n\equiv (\sum_{r=1}^{n} b_{r,j} ) \mod 2$, then $b_{n+1,j}:=1$.\\
Construct rows $n+2$ to $2(n+1)$: \\ \phantom{x} for $i$ from $1$ to $n+1$ do: \\ \phantom{xx} for $j$ from $1$ to $2^n$ do: \phantom{xxx}
$b_{n+1+i},j:=1-b_{i,j}$.
\end{alg}

One checks that this algorithm gives indeed the matrices $B_n$ as defined in Section 3.

The $(n+1+i)^{th}$ row $r_{n+1+i}$ of $B_n$ is by definition the binary complement of the $i^{th}$ row $r_i$ of $B_n$.  Suppose that $r_i \cdot k=0$ for some
vector $k$.  Since all entries of $B_n$ are nonnegative, a subvector of $k$ restricted to the entries where $r_i$ is nonzero must be homogeneous in the sense
that the sum of the positive entries equals the sum of the negative entries. But since the ideal $I_{L_n}$ itself is homogeneous (\cite{St}), the same must be
true for the subvector of $k$ restricted to the entries where $r_i$ is zero. Hence $r_{n+1+i}\cdot k=0$.  Therefore, it is enough to analyze the top half of the
matrix $B_n$ when determining the kernel elements.

\begin{remark}
There are $n$ copies of $B_{n-1}$ inside $B_n$.\\
By deleting one leaf at a time, we get $n$ copies of $T_{n-1}$ as a subtree of $T_n$. Suppose we delete leaf
$(i)$ from $T_n$ to get the tree $T_n^{(i)}$ on leaves $1,2,\dots,i-1,i+1,\dots,n$. Ignoring the two rows of $B_n$ that represent the leaf $(i)$ and taking into
account the columns of $B_n$ containing nonzero entries of the row  indexed by $a_0^{(i)}$ (that is, observing $0$ at leaf $(i)$) gives precisely the matrix
$B_{n-1}$ corresponding to $T_n^{(i)}$. Note that the entry indexed by $a^{(n+1)}_{g}$, for any $g\in G$, will be correct since we are ignoring the identity of
the group, as in Remark \ref{obs:1}.
\end{remark}

This leads to a way of constructing  a basis of $\ker (B_n)$ from the one of $\ker(B_{n-1})$. Namely,  removing leaf $(1)$ from $T_n$ produces
$\dim(\ker(B_{n-1}))=2^{n-1}-n-1$ independent vectors in $\ker(B_n)$. Let us name this collection of vectors $V_1$.  Removing leaf $(2)$ produces a collection
$V_2$ consisting of $\dim(\ker{B_{n-1})} - \dim(\ker{B_{n-2})}= 2^{n-2}-1$ vectors in $\ker(B_n)$.  $V_2$ is independent of $V_1$ since the second half of each
vector in $V_2$ has nonzero entries in the columns of $B_n$ where all vectors in $V_1$ are zero, a direct consequence of the location of the submatrix
corresponding to $T_n^{(2)}$.  Finally,  removing any other leaf $(i)$ of $T_n$ produces a collection $V_i$ of as many new kernel elements as there are new
columns involved (in terms of the submatrix structure); namely, $2^{n-i}$ new vectors.  Note that every vector in $V_2$ has a nonzero entry in at least one new
column so that the full collection is independent of $V_1$.

Using the above procedure, we have obtained $$ (2^{n-1}-n-1) + (2^{n-2}-1) + (2^{n-3}) + \dots + 2^{n-n} $$ independent vectors in the kernel of $B_n$.  This is
exactly one less than the desired number, $2^n-n-2$.  Hence to the list of the kernel generators we add one additional vector $v$ that is independent of all the
$V_i$, $i=1,\dots,n$ as it has a nonnegative entry in the last column.  (Note that no $v\in V_i$ has this property by the observation on the column location of
the submatrix associated with each $T_n^{(i)}$.) In particular,  $v= [0,\dots,0 ,1,0,0,-1,-1,0,0,1]\in \ker (B_n).$  To see this, we simply notice that the rows
of  the last $8$-column block of $B_n$ are precisely the rows of the first $8$-column block of $B_n$ up to permutation of rows, which does not affect the kernel.

The lattice basis we just constructed is directly computed by the following algorithm.

\begin{alg}\rm  \label{alg1}
 [Construction of the lattice basis for $T_n$]\\
{\bf Input:} the number of leaves $n$ of the claw tree $T_n$.\\ {\bf Output:} a basis of $\ker B_n$ in form of a $(2^n-n-2)\times 2^n$ matrix $L_n$.
\\
Let $L_3:=\left[\begin{matrix} 0&0&1&-1&-1&1&0&0\\ 0&1&0&-1&-1&0&1&0\\1&0&0&-1&-1&0&0&1\end{matrix}\right]$.\\ Set $k:=4$. \\ The following subroutine lifts
$L_{k-1}$ to $L_k$: \\ WHILE $k\leq n$ do:\{
\\ Initialize $L_k$ to the zero matrix.\\
 For $i$ from $1$ to $k$ do: \\ \phantom{x} $\cols (i):=\{1.. 2^{k-i},(2)2^{k-i}+1.. (3)2^{k-i},\dots, (2^i-2)2^{k-i}+1..
(2^i-1)2^{k-i} \}$.\\  Denote by $L_{k,j}[\cols(i)]$ the $j^{th}$ row vector of $L_k$ restricted to columns $\cols(i)$. \\ Set $i:=1$: \\ \phantom{x} for $j$
from $1$ to $2^{k-1}-k-1$ do:  \phantom{xxx} $L_{k,j}[\cols(i)]:=L_{k-1,j}$.\\ Set $i:=2$: \\ \phantom{x} for $j$ from $1$ to $2^{k-2}-1$ do :
\\ \phantom{xx}$L_{k,(2^{k-1}-k-1) + j}[\cols(i)] := L_{k-1,(2^{k-1}-k-1)-(2^{k-2}-1)+j}$.\\  For $i$ from $3$ to $k$ do: \\ \phantom{x} for $j$ from $1$ to
$2^{k-i}$ do: \\ \phantom{xx} $L_{k,(2^k-2^{k+1-i}-k-2)+j}[\cols(i)]:=L_{k-1, (2^{k-1}-k-1)-(2^{k-i})+j}$.\\ Finally,
$L_{k,2^k-k-2}[2^k-7..2^k]:=[1,0,0,-1,-1,0,0,1]$.\\ RETURN $L_k$. \}
\end{alg}

\begin{example}
Consider the tree on $n=4$ leaves. Then
 $$B_4 =
 \left[ \begin {array}{cccccccccccccccc} 1&1&1&1&1&1&1&1&0&0&0&0&0&0&0
&0\\\noalign{\medskip}1&1&1&1&0&0&0&0&1&1&1&1&0&0&0&0
\\\noalign{\medskip}1&1&0&0&1&1&0&0&1&1&0&0&1&1&0&0
\\\noalign{\medskip}1&0&1&0&1&0&1&0&1&0&1&0&1&0&1&0
\\\noalign{\medskip}1&0&0&1&0&1&1&0&0&1&1&0&1&0&0&1
\\\noalign{\medskip}0&0&0&0&0&0&0&0&1&1&1&1&1&1&1&1
\\\noalign{\medskip}0&0&0&0&1&1&1&1&0&0&0&0&1&1&1&1
\\\noalign{\medskip}0&0&1&1&0&0&1&1&0&0&1&1&0&0&1&1
\\\noalign{\medskip}0&1&0&1&0&1&0&1&0&1&0&1&0&1&0&1
\\\noalign{\medskip}0&1&1&0&1&0&0&1&1&0&0&1&0&1&1&0\end {array}
 \right].$$
The lattice basis is given by the rows of the following matrix:

$$L_4=\left[
\begin {array}{cccccccccccccccc}0 &0 &1 &-1 &-1 & 1 & 0 &0  & 0 &0 &0 & 0& 0 &0 &0 &0
\\\noalign{\medskip}0 &1 &0 &-1 &-1 & 0 & 1 &0  & 0 &0 &0 & 0& 0 &0 &0 &0
\\\noalign{\medskip}1 &0 &0 &-1 &-1 & 0 & 0 &1  & 0 &0 &0 & 0& 0 &0 &0 &0
\\\noalign{\medskip}0 &0 &1 &-1 & 0 & 0 & 0 &0  &-1 &1 &0 & 0& 0 &0 &0 &0
\\\noalign{\medskip}0 &1 &0 &-1 & 0 & 0 & 0 &0  &-1 &0 &1 & 0& 0 &0 &0 &0
\\\noalign{\medskip}1 &0 &0 &-1 & 0 & 0 & 0 &0  &-1 &0 &0 & 1& 0 &0 &0 &0
\\\noalign{\medskip}0 &1 &0 & 0 & 0 &-1 & 0 &0  &-1 &0 &0 & 0& 1 &0 &0 &0
\\\noalign{\medskip}1 &0 &0 & 0 & 0 &-1 & 0 &0  &-1 &0 &0 & 0& 0 &1 &0 &0
\\\noalign{\medskip}1 &0 &0 & 0 & 0 & 0 &-1 &0  &-1 &0 &0 & 0& 0 &0 &1 &0
\\\noalign{\medskip}0 &0 &0 & 0 & 0 & 0 & 0 &0  & 1 &0 &0 &-1&-1 &0 &0 &1\end {array}
\right]. $$ The lattice vectors correspond to the relations on the leaf observations in the natural way; namely, the first column corresponds to $q_{0,\dots,0}$,
the second to $q_{0,\dots,0,1}$, and so on. Therefore, the lattice basis ideal for $T_4$ in Fourier coordinates is

\begin{align*}
I_{L_4} = (
&q_{0010}q_{0101}-q_{0011}q_{0100},
q_{0001}q_{0110}-q_{0011}q_{0100},
q_{0000}q_{0111}-q_{0011}q_{0100},\\
&q_{0010}q_{1001}-q_{0011}q_{1000},
q_{0001}q_{1010}-q_{0011}q_{1000},
q_{0000}q_{1011}-q_{0011}q_{1000}, \\
&q_{0001}q_{1100}-q_{0101}q_{1000},
q_{0000}q_{1101}-q_{0101}q_{1000},  \\
&q_{0000}q_{1110}-q_{0110}q_{1000},
q_{1000}q_{1111}-q_{1011}q_{1100} ).
\end{align*}
\end{example}
This ideal is contained in the ideal of phylogenetic invariants $I_4$ for $T_4$.  In the next section, we compute explicitly the generators of the ideal of
invariants for any claw three $T_n$ and the group $\mathbb Z_2$.

\section{Ideal of invariants}

We show that the lattice basis ideals provide basic building blocks for the full ideals of invariants, as expected. However, instead of computing the ideal of
invariants as a saturation of the lattice basis ideal in a standard way (e.g. \cite{HeMa},\cite{St}), we use the recursive constructions from the previous
section on the saturated ideals directly.  We begin with the ideal of invariants for the smallest tree, and build all other trees recursively. The underlying
ideas for how to lift the generating sets come from Algorithm \ref{alg1}.\\ We will denote the ideal of the claw tree on $n$ leaves by $I_n=\ker\varphi_n$. As we
have seen, the first nontrivial ideal is $I_3$.

\subsection{The tree on $n=3$ leaves}
\begin{claim}
The ideal of the claw tree on $n=3$ leaves is
$$I_3=(q_{000}q_{111}-q_{100}q_{011},q_{001}q_{110}-q_{100}q_{011},q_{010}q_{101}-q_{100}q_{011}).$$
\end{claim}

This can be verified by computation. In particular, this ideal is equal to the lattice basis ideal for the tree on three leaves; $I_{L_3}$ is already prime in
this case.

Let $<:=<_{lex}$ be the lexicographic order on the variables induced by $$q_{000}>q_{001}>q_{010}>q_{011}>q_{100}>q_{101}>q_{110}>q_{111}.$$  (That is,
$q_{ijk}>q_{i'j'k'}$ if and only if $(ijk)_2<(i'j'k')_2$, where $(ijk)_2$ denotes the binary number $ijk$.)
\begin{remark}
The three generators of $I_3$ above are a Gr\"{o}bner basis for $I_3$ with respect to $<$, since the initial terms, written with coefficient $+1$ in the above
description, are relatively prime so all the S-paris reduce to zero.
\end{remark}

\begin{remark} Write the quadratic binomial $q=q^+ - q^- $ as
$$ q_{g_1^{(1)}g_1^{(2)}g_1^{(3)}}q_{g_2^{(1)}g_2^{(2)}g_2^{(3)}}-q_{h_1^{(1)}h_1^{(2)}h_1^{(3)}}q_{h_2^{(1)}h_2^{(2)}h_2^{(3)}}.$$
Then $q \in I_3 $ if and only if the following two conditions hold:
\begin{enumerate}
\item  Exchanging the roles of $q_{h_1^{(1)}h_1^{(2)}h_1^{(3)}}$ and $q_{h_2^{(1)}h_2^{(2)}h_2^{(3)}}$ if necessary,
$$g_1^{(1)}+g_1^{(2)}+g_1^{(3)}=h_1^{(1)}+h_1^{(2)}+h_1^{(3)} $$ and
$$g_2^{(1)}+g_2^{(2)}+g_2^{(3)}=h_2^{(1)}+h_2^{(2)}+h_2^{(3)}, $$
\item $g_1^{(i)}+g_2^{(i)}=1=h_1^{(i)}+h_2^{(i)}$ for $1\leq i\leq 3=n$.
\end{enumerate}
\end{remark}

Note that the second condition holds since otherwise the projection of $q$ obtained by eliminating the leaf $(i)$ at which the observations $g_1^{(i)}$ and
$g_2^{(i)}$ are both equal to $0$ or to $1$ produces an element $q'$ in the kernel of the map $\varphi_2$ of the $2$-leaf tree, which is trivial.

\subsection{The tree on an arbitrary number of leaves}
\paragraph{}Let us now define a set of maps and a distinguished set of binomials in $I_n$.
\begin{definition}
Let  $\pi_i(q)$ be the projection of $q$ that eliminates the $i^{th}$ index of each variable in $q$.
\end{definition}
For example, $$\pi_4(q_{0000}q_{1110}-q_{1000}q_{0110})=q_{000}q_{111}-q_{100}q_{011} .$$

\begin{definition}  Assume that $n\geq 4$.

Let $\mathcal G_n$ be the set of quadratic binomials $q\in I_n$ that can be written as
$$q=q^+ - q^- =  q_{g_1^{(1)}\dots g_1^{(n)}}q_{g_2^{(1)}\dots g_2^{(n)}} -  q_{h_1^{(1)}\dots h_1^{(n)}} q_{h_2^{(1)}\dots h_2^{(n)}}$$
 such that one of the two following properties is satisfied:
\newcounter{Lcount}
  \begin{list}{Property (\roman{Lcount}):}
    {\usecounter{Lcount}
    \setlength{\rightmargin}{\leftmargin}}
\item
  For {\em some} $1\leq i\leq n$, $j\in\mathbb Z_2$,
  \begin{equation}\label{eq:i}
      g_1^{(i)}=g_2^{(i)}= j =h_1^{(i)}=h_2^{(i)}
      \end{equation}
 and
  \begin{equation}\label{eq:ia}
   \pi_i(q)\in I_{n-1}.
   \end{equation}
\item
 For {\em each} $1\leq k\leq n$,
       \begin{equation}\label{eq:ii}
       g_1^{(k)}+g_2^{(k)} = 1 = h_1^{(k)}+h_2^{(k)}
       \end{equation}
and \begin{equation}\label{eq:iia}
   \pi_k(q)\in I_{n-1}.
    \end{equation}
\end{list}
\end{definition}

\begin{example}
Let $n=4$.  The set of elements $q\in\mathcal G_n$  with Property (i) consists of those for which $j=0$:\\
$q_{0000}q_{0111}-q_{0100}q_{0011}$, $q_{0001}q_{0110}-q_{0100}q_{0011}$, $q_{0010}q_{0101}-q_{0100}q_{0011}$,\\
$q_{0000}q_{1011}-q_{1000}q_{0011}$, $q_{0001}q_{1010}-q_{1000}q_{0011}$, $q_{0010}q_{1001}-q_{1000}q_{0011}$,\\
$q_{0000}q_{1101}-q_{1000}q_{0101}$, $q_{0001}q_{1100}-q_{1000}q_{0101}$, $q_{0100}q_{1001}-q_{1000}q_{0101}$,\\
$q_{0000}q_{1110}-q_{1000}q_{0110}$, $q_{0010}q_{1100}-q_{1000}q_{0110}$, $q_{0100}q_{1010}-q_{1000}q_{0110}$;\\
and those for which $j=1$:\\
$q_{1000}q_{1111}-q_{1100}q_{1011}$, $q_{1001}q_{1110}-q_{1100}q_{1011}$, $q_{1010}q_{1101}-q_{1100}q_{1011}$,\\
$q_{0100}q_{1111}-q_{1100}q_{0111}$, $q_{0101}q_{1110}-q_{1100}q_{0111}$, $q_{0110}q_{1101}-q_{1100}q_{0111}$,\\
$q_{0010}q_{1111}-q_{1010}q_{0111}$, $q_{0011}q_{1110}-q_{1010}q_{0111}$, $q_{0110}q_{1011}-q_{1010}q_{0111}$,\\
$q_{0001}q_{1111}-q_{1001}q_{0111}$, $q_{0011}q_{1101}-q_{1001}q_{0111}$, $q_{0101}q_{1011}-q_{1001}q_{0111}$.\\
 The set of elements $q\in\mathcal G_n$ with Property (ii) are: \\
$q_{0000}q_{1111}-q_{1001}q_{0110}$, $q_{0001}q_{1110}-q_{1000}q_{0111}$, $q_{0011}q_{1100}-q_{1001}q_{0110}$,\\
$q_{0010}q_{1101}-q_{1000}q_{0111}$, $q_{0101}q_{1010}-q_{1001}q_{0110}$, $q_{0100}q_{1011}-q_{1000}q_{0111}$.
\end{example}

\begin{proposition}\label{gens}
For $n\geq 4$, the set of binomials in $\mathcal G_n$ generates the ideal $I_n$. That is,
$$ I_n = (q:q^+-q^-\in\mathcal G_n).$$
In addition, this set of generators can be obtained inductively by lifting the generators corresponding to the various phylogenetic ideals on $n-1$ leaves.
\end{proposition}
\begin{proof}
Condition (\ref{eq:ii}) is simply the negation of (\ref{eq:i}). Condition (\ref{eq:i}) can be restated as follows: for some $1\leq i\leq n$ and a fixed $j$,
$$(a_j^{(i)})^2|\varphi_n(q^+)\mbox{   and   } (a_j^{(i)})^2|\varphi_n(q^-). $$ Therefore, Property (i) translates to having an observation $j$ fixed at leaf
$(i)$ for each of the variables in $q$.
On the other hand, condition (\ref{eq:ii}) means that for any $k$, not all the $k^{th}$ indices are $0$ and not all are
$1$.  Thus Property (ii) means that no leaf has a fixed observation, and can be restated as follows: for every $1\leq i\leq n$,
 \begin{equation}\label{eq:ii*}
a_0^{(i)}a_1^{(i)}|\varphi_n(q^+) \mbox{   and   } a_0^{(i)}a_1^{(i)}|\varphi_n(q^-).
\end{equation}

By definition, the ideal $I_n$ is toric, so it is generated by binomials. In fact, it is generated by homogeneous binomials, because each row of the matrix $B_n$
used for defining it has row sum $n+1$ (\cite{St}, chapter $4$).  In addition, Sturmfels and Sullivant in \cite{StSu} have shown that the ideal $I_n$ is
generated in degree $2$.  Hence it suffices to consider homogeneous quadratic binomials.  Let $q=q^+-q^-$ be a binomial in $I_n$ of degree $2$. Then clearly
either (\ref{eq:i}) or (\ref{eq:ii}) holds; that is, either the index corresponding to one leaf is fixed for all the monomials in $q$, or none of them are. \\

 In the former case, for the index $i$ from equation (\ref{eq:i}),
\begin{align*}
q\in I_n & \iff \varphi_n(q^+)=\varphi_n(q^-)\\
 & \iff \varphi_{n-1}(\pi_i(q^+)) = \varphi_{n-1}(\pi_i(q^-)) \iff \pi_i(q)\in I_{n-1},
\end{align*}
where the first statement holds by definition of $\varphi_n$ and the second by definition of the projection $\pi_i$.\\

In the latter case, for each $i$ with $1\leq i \leq n$,
\begin{align*}
q\in I_n & \iff \varphi_n(q^+)=\varphi_n(q^-) \\
 & \iff \varphi_{n-1}(\pi_i(q^+))=\varphi_{n-1}(\pi_i(q^-))\iff \pi_i(q)\in I_{n-1},
\end{align*}
 where the second statement holds by definition of $\pi_i$ and (\ref{eq:ii*}). It follows that $I_n=(q: q\in\mathcal G_n)$.\\

In particular, the set of generators for $I_n$ with Property (i) can be obtained from those of $I_{n-1}$ by inserting first $0$ at the $i^{th}$ index position
for each monomial of $q\in\mathcal G_{n-1}$ and then repeating the same process by inserting $1$. This operation corresponds to lifting to all the possible
preimages of $\pi_i(q)$ that satisfy Property (i) for each $1\leq i\leq n$ and every $q\in\mathcal G_{n-1}$.  The set of generators for $I_n$ with Property (ii)
can be obtained from those of $I_{n-1}$ by a  similar lifting to all preimages of $\pi_i(q)$ for each $q\in\mathcal G_{n-1}$ in such a way that Property (ii) is
satisfied.  Namely, for every $q=q^+-q^-\in\mathcal G_{n-1}$ with Property (ii), one inserts $0$ at the $i^{th}$ index position for one monomial of $q^+$ and for
one monomial of $q^-$, and inserts $1$ at the $i^{th}$ index position for the remaining monomials of $q^+$ and $q^-$. In addition, by definition of Property
(ii), it suffices to lift to the preimages of $\pi_n(q)$ only. \hfill $\square$
\end{proof}

\begin{remark}
A different recursion has been proposed by Sturmfels and Sullivant in \cite{StSu2}.
\end{remark}

\paragraph{}

Recall (\cite{St}) that a binomial $q=q^+-q^-\in I$ is said to be {\em primitive} if there exists no binomial $f=f^+-f^-\in I$ with the property that $f^+|q^+$
and $f^-|q^-$. A {\em circuit} is a primitive binomial of minimal support.
\begin{remark}
The binomials in $\mathcal G_n$ are circuits of $I_n$,  since the ideal is generated by squarefree binomials and contains no linear forms.
\end{remark}

In general, we can describe the generators of $I_n$ as follows: given $n$, begin by lifting $\mathcal G_3$ recursively  to produce $\mathcal G_{n-1}$; that is,
until the number of indices of each generator reaches $n-1$. Next, lift $\mathcal G_{n-1}$ $n$ times so that Property (i) is satisfied for one of the $n$ index
positions. For example, $$q:=q_{0000}q_{1111}-q_{1001}q_{0110}\in\mathcal G_4$$ can be lifted to a generator of $I_5$ in ten different ways: by lifting to
preimages of $\pi_1,\dots, \pi_5$ so that Property (i) is satisfied with either a $0$ or a $1$:
\begin{align*}
\pi_1^{-1}(q)& = \{ q_{00000}q_{01111}-q_{01001}q_{00110}, q_{10000}q_{11111}-q_{11001}q_{10110} \},\\
\pi_2^{-1}(q)& = \{ q_{00000}q_{10111}-q_{10001}q_{00110} , q_{01000}q_{11111}-q_{11001}q_{01110} \},
\end{align*}
and so on.
This will be the set of binomials in $\mathcal G_n$ with Property (i).  Clearly, some generators will repeat during the recursive lifting: lifting by
inserting $0$ at position $(i)$ allows the $0$ to occur at the previous $i-1$ positions. Also, fixing $1$ at any leaf allows $0$ to appear on any of the other
leaves.

To construct $q^+-q^-$ with Property (ii), we need not proceed inductively, as all projections of binomials that satisfy this property must satisfy it, too.
Instead, we consider two cases corresponding to the parity of $n$. Namely, recalling the definition of Property (ii), first we fix $q^-$ in such a way to  ensure
that $in_{<_{lex}}(q)=q^+$.

Suppose $n$ is odd.  Fix $q^-$ by taking $$q^- = q_{01\dots 1}q_{10\dots 0}$$ with $n$ indices in each of the two variables. Then $n-1$ being even provides that
$a_0^{(n+1)}a_1^{(n+1)} | \varphi_n (q^-)$. Thus every choice of $q^+$ must satisfy the same.  To find $q^+$, we need to choose pairs of $n$-digit binary numbers
with digits complementary to each other, and thus there are $2^{n-1}-1$ choices for $q^+$.  Specifically, listing the smallest $2^{n-1}-1$ $n$-digit binary
numbers and pairing them with the largest $2^{n-1}-1$ $n$-digit binary numbers in reverse order produces all choices for $q^+$, and we have a complete list of
generators.  For example, the first such generator in the list would be $q_{0\dots 0}q_{1\dots 1}-q_{01\dots 1}q_{10\dots 0}$.

If $n$ is even, then we can create $q^-$ such that $(a_0^{(n+1)})^2$ or $(a_1^{(n+1)})^2$ divides $\varphi_n(q^-)$ and $\varphi_n(q^+)$. Namely, the two choices
for $q^-$ are $$q^-=q_{01\dots 1}q_{10\dots 0} \mbox{  and  } q^-=q_{01\dots 10}q_{10\dots 01}.$$ The list of all possible $q^+$ is obtained in the manner
similar to the case when $n$ is odd, except that the odd pairs in the list receive the first choice of $q^-$, while the even pairs receive the second. The number
of such generators $q^+-q^-$ is $2^{n-1}-2$, since there are $2^n$ $n$-digit binary numbers and thus half as many pairs, and $2$ choices are taken by the $q^-$.

In summary, the number of generators of $I_n$ that satisfy Property (ii) is \\
$(2^{n-1}-2) +(n\mod 2)$.

Next we strengthen Proposition (\ref{gens}).

\begin{proposition}\label{groebner}
The  set $\mathcal G_n$ is a lexicographic Gr\"{o}bner basis of $I_n$, for any $n\geq 4$.
\end{proposition}
\begin{proof}

For the case $n=3$ this is already shown. Let $n>3$. Then we can partition the set of $q\in\mathcal G_n$ into those satisfying Property (i) or (ii). Note that
$I_n$ is prime by definition, and thus radical. Also, Proposition (\ref{gens}) shows it is generated by squarefree quadratic binomials.  These facts are used in
what follows.\\

Let $q_i$,$q_j\in I_n$.  If $(q_i^+,q_j^+)=1$, the S-pair $S(q_i,q_j)$ reduces to zero. Also, if $q_i^-$ and $q_j^-$ are not relatively prime, the cancellation
criterion provides that the corresponding S-pair also reduces to zero. Therefore we consider $f:=S(q_i,q_j)\in I_n$ with $(q_i^+,q_j^+)\neq 1$ and
$(q_i^-,q_j^-)=1$. In particular, $\deg(f)=3$.   Let us write $q_i=q_{g_1}q_{g_2}-q_{h_1}q_{h_2}$ and $q_j=q_{g_1}q_{g_3}-q_{h_3}q_{h_4}$. Then
$$f=q_{g_3}q_{h_1}q_{h_2}-q_{g_2}q_{h_3}q_{h_4}\in I_n.$$ {\em Case I.} Suppose $q_i$ satisfies Property (i) and $q_j$ satisfies Property (ii). Then there exists
a $k$ such that $\pi_k(q_i)\in I_{n-1}$. Furthermore, Property (ii) implies that $\pi_k(q_j)\in I_{n-1}$. A very technical argument  shows that $$\pi_k(f)\in
I_{n-1}$$ and furthermore, this projection preserves the initial terms. In summary, to check that $\pi_k(f)\in I_{n-1}$, it suffices to ensure that
$a_s^{(n)}|\varphi_{n-1}(\pi_k(q_{g_3}q_{h_1}q_{h_2}))$ if and only if $a_s^{(n)}|\varphi_{n-1}(\pi_k(q_{g_2}q_{h_3}q_{h_4}))$, where $s$ is the sum of the
observations on the leaves of the $(n-1)$-leaf tree obtained from $T$ by deleting leaf $(k)$.  There are two cases corresponding to the parity of $n$. If $n$ is
odd, there are additional subcases determined by the correspondence  of the images of the variables in the two monomials of $f$ under $\varphi_{n-1}$. The facts
that $q_i$ and $q_j$ satisfy Properties (i) and (ii), respectively, play a crucial role in the argument.  Checking all the cases then shows that $\pi_k(f)\in
I_{n-1}$ and that initial terms are preserved under this projection.

 Applying the induction hypothesis then finishes the proof.\\
{\em Case II.} Suppose both $q_i$ and $q_j$ satisfy Property (i).  Then there is a $q_k\in\mathcal G_n$ satisfying Property (ii) where both $S(q_i,q_k)$ and
$S(q_j,q_k)$ reduce to zero.  The three-pair criterion (\cite{HeMa}) provides the desired result.\\ {\em Case III.} If both $q_i$ and $q_j$ satisfy Property
(ii), then it can be seen from the construction preceding this Proposition that the initial terms are relatively prime, so their S-polynomial need not be
considered. \hfill$\square$ \end{proof}

Proposition \ref{groebner} has important theoretical consequences.  Let $S$ be a polynomial ring over the field $K$. Recall (\cite{CRV}) that $S/I$ is {\em
Koszul} if the field $K$ has a linear resolution as a graded $S/I$-module: $$\dots\rightarrow (S/I)^{\beta_2}(-2)\rightarrow (S/I)^{\beta_1}(-1)\rightarrow
S/I\rightarrow K\rightarrow 0.$$  An ideal $I\subset S$ is said to be quadratic if it is generated by quadrics.  $S/I$ is {\em quadratic} if its defining ideal
$I$ is quadratic, and it is {\em G-quadratic} if $I$ has a quadratic Gr\"{o}bner basis.  It is known (e.g. \cite{CRV}) that if $S/I$ is  G-quadratic, then it is
Koszul, which in turn implies it is quadratic.  The reverse implications do not hold in general.
 We have just found an infinite family of toric varieties whose coordinate rings $S/I$ are G-quadratic.
\begin{corollary}
The coordinate ring of the toric variety whose defining ideal is $I_n$ is Koszul for every $n$.
\end{corollary}

The approach developed here produces the list of generators for the kernel of $B_n$ all of which are of degree two.  In addition, by constructing the toric
ideals of invariants inductively, we are able to explicitly calculate the quadratic Gr\"{o}bner bases.  In light of the conjecture posed in \cite{StSu} that the
ideal of phylogenetic invariants for the group of order $k$ is generated in degree at most $k$, we are working on generalizing the above approach to any abelian
group of order $k$. In particular, we want to give a description of  the lattice basis ideal $I_{L_n}$ and the ideal of invariants $I$ for $G=\mathbb Z_2 \times
\mathbb Z_2$ with generators of degree at most $4$. These phylogenetic ideals are of interest to computational biologists.


\begin{thebibliography}{12}
\bibitem{AlRh}
E. Allman, J. Rhodes: Phylogenetic ideals and varieties for the general Markov model. Advances in Applied Mathematics, to appear. Preprint,
arXiv.org:math/0410604, 2004.

\bibitem{AlRh2}
E. Allman, J. Rhodes: Identifying evolutionary trees and substitution parameters
        for the general Markov model with invariable sites.  Preprint, arXiv.org:q-bio/0702050, 2007.

\bibitem{CaFS}
M. Casanellas, J. Fernandez-Sanchez:  Performance of a new invariants method on homogeneous and non-homogeneous quartet trees. Preprint, arXiv.org:q-bio/0610030,
2006.

\bibitem{CRV}
A. Conca, M.E. Rossi, G. Valla: Gr\"{o}bner flags and Gorenstein algebras. Compositio Math. {\bf 129} (27), Number 1, October 2001, pp. 95--121

\bibitem{Eis}
D. Eisenbud: Commutative algebra with a view toward algebraic geometry. Graduate Texts in Mathematics 150, Springer-Verlag, 1995.

\bibitem{ErSt}
N. Eriksson, K. Ranestad, B. Sturmfels, S. Sullivant: Phylogenetic Algebraic Geometry. In:  Projective varieties with unexpected properties, (editors C.
Ciliberto, A. Geramita, B. Harbourne, R-M. Roig and K. Ranestad), De Gruyter, Berlin, 2005, pp. 237-255

\bibitem{ErYa}
N. Eriksson, Y. Yao: Metric learning for phylogenetic invariants. Preprint, arXiv.org:q-bio/0703034, 2007.

\bibitem{HeMa}
R. Hemmecke, P. Malkin:
Computing generating sets of lattice ideals.
Preprint, arXiv.org:math/0508359, 2005

\bibitem{PaSt}
L. Pachter, B. Sturmfels: Algebraic statistics for computational biology. Cambridge University Press, New York, NY, USA, 2005

\bibitem{St}
 B. Sturmfels:
 Gr\"{o}bner bases and convex polytopes.
 American Mathematical Society, University Lecture Series {\bf 8}, 1996

\bibitem{StSu}
B. Sturmfels, S. Sullivant: Toric ideals of phylogenetic invariants. J. Comp. Biol. {\bf 12} (2005), pp. 204-228.

\bibitem{StSu2}
B. Sturmfels, S. Sullivant:
Toric geometry of cuts and splits.
Preprint, arXiv.org:math.AC/0606683, 2006.

\end{thebibliography}
\end{document}